\documentclass[conference]{IEEEtran}
\IEEEoverridecommandlockouts
\usepackage{cite}
\usepackage{amsmath,amssymb,amsfonts}
\usepackage{algorithmic}
\usepackage{booktabs}
\usepackage[pdftex]{graphicx}
\usepackage{textcomp}
\usepackage[table]{xcolor}
\usepackage{subcaption}
\usepackage{multirow}
\def\BibTeX{{\rm B\kern-.05em{\sc i\kern-.025em b}\kern-.08em
    T\kern-.1667em\lower.7ex\hbox{E}\kern-.125emX}}

\graphicspath{{./images/}}

\definecolor{Gray}{gray}{0.95}
\newcolumntype{g}{>{\columncolor{Gray}}c}
\begin{document}

\title{Numerical computing in engineering mathematics
\thanks{© 2022 IEEE.  Personal use of this material is permitted.  Permission from IEEE must be obtained for all other uses, in any current or future media, including reprinting/republishing this material for advertising or promotional purposes, creating new collective works, for resale or redistribution to servers or lists, or reuse of any copyrighted component of this work in other works}
}

\author{\IEEEauthorblockN{Firuz Kamalov}
\IEEEauthorblockA{\textit{Department of Electrical Engineering} \\
\textit{Canadian University Dubai}\\
Dubai, UAE \\
firuz@cud.ac.ae}
%\and
%\IEEEauthorblockN{Author2}
%\IEEEauthorblockA{\textit{Department of Author2} \\
%\textit{University of Author2}\\
%Author2City, UAE \\
%Author2@email.com}
\and
\IEEEauthorblockN{Ho-Hon Leung}
\IEEEauthorblockA{\textit{Department of Mathematical Sciences} \\
\textit{United Arab Emirates University}\\
Al Ain, UAE \\
hohon.leung@uaeu.ac.ae}
}

\maketitle

\begin{abstract}
The rapid advances in technology over the last decade have significantly altered the nature of engineering knowledge and skills required in the modern industries. In response to the changing professional requirements, engineering institutions have updated their curriculum and pedagogical practices. However, most of the changes in the curriculum have been focused on the core engineering courses without much consideration for the auxiliary courses in mathematics and sciences. In this paper, we aim to propose a new, augmented mathematics curriculum aimed at meeting the requirements of the modern, technology-based engineering workplace. The proposed  updates require minimal resources and can be seamlessly integrated into the existing curriculum.
\end{abstract}


%----------------------------------------------------------------------------------------------------------------------------------
%----------------------------------------------------------------------------------------------------------------------------------
\section{Introduction}
The 4th Industrial Revolution has had a dramatic impact on the engineering profession. The modern technologies such as artificial intelligence, the internet of things, and advanced robotics have altered engineering systems and processes. Today's engineers are expected to be able to leverage these resources to produce their products. To meet the new professional requirements, engineering educational institutions have revised their curricula. The changes in the curricula include both updating the existing programs as well as introducing completely new programs. Given the rapid technological progress, universities and colleges around the world are continuously adapting to the ever-changing environment. While a significant progress in modernizing the engineering curriculum has been achieved, there still remains room for improvement. 

Catalyzed by  the exponential increase in computational power and interconnectedness, the modern industrial revolution has reshaped the skills and competencies required of the engineers. The changes in engineering curricula in response to Industry 4.0 have been threefold: i) modernizing the existing programs, ii) introduction of new programs, and iii) revising the pedagogical approach. Modernizing the existing programs involves introduction of new courses in the study plan related to emerging technologies. In addition, existing courses can be updated with new content. Fresh new programs in emerging technologies are also introduced by universities and colleges. Many institutions now offer degrees in artificial intelligence and mechatronics which were not there 20 years ago. Finally, universities have revised their approaches to course delivery. Student-centered learning, project-based learning, and applied learning have become popular in the new engineering educational paradigm.

While significant effort has been made to revise the core engineering courses, the auxiliary courses in mathematics and sciences received little consideration.
The mathematics and sciences courses play a key role in the engineering curriculum. Given their importance, the curriculum updates must also be extended to the auxiliary courses. By implementing a comprehensive update of the engineering curriculum that includes both the core and auxiliary courses, a more effective outcome can be achieved. 

The goal of this paper is to propose a modernized engineering mathematics curriculum in line with the broader efforts to update engineering education to adapt to Industry 4.0. The key feature of the new curriculum is the introduction numerical computing in the existing mathematics courses. The latest industrial revolution has been driven largely  by the dramatic increase in computational power. Therefore, today's engineers must be well-equipped to leverage the computing power in their work. 

Since mathematics courses are usually taken at the beginning of the study plan, it offers a natural avenue for introducing numerical computing to students. Furthermore, many problems in mathematics can be solved numerically making it natural to integrate numerical computing in mathematics courses. By studying numerical computing in mathematics courses, students will acquire the necessary theoretical and practical skills to apply in their downstream, specialized engineering courses.

This paper is structured as follows. Section 2 provides an overview of the existing efforts to update the engineering curriculum in response to Industry 4.0. Section 3 discusses the current approaches to integrate scientific computing in mathematics courses. In Section 4, we present our proposal for modernizing the mathematics curriculum to integrate numerical computing. Section 5 concludes the paper with final remarks.

\section{Engineering education and Industry 4.0}
Engineering departments in colleges and universities have made significant changes in their curricula in response to the new environment created by the recent, rapid advances in technology.
In particular, the existing programs have been updated to include courses that target  emerging technologies. Completely new programs related to AI and mechatronics have also been adopted by universities. Innovations in the field of engineering education continue to take place with new developments on the horizon.

There exist several studies investigating the modern engineering curricula and evaluating their effectiveness. It is argued in \cite{Hadgraft} that engineering educators must prepare their students to face three key challenges: sustainability, the 4th Industrial Revolution, and employability. The authors find that colleges and universities are responding to these challenges by emphasizing student-centered learning, integration of theory and practice, digital and online learning, and the definition of professional competencies. In particular, response to the needs of Industry 4.0 require interdisciplinary collaboration across several programs and disciplines. Interaction and integration of technologies plays a key role in this process \cite{Kamalov2, Lorenz}. Interdisciplinary engineering education  requires sound pedagogy and teaming experiences to encourage student in collaborative and interdisciplinary practice \cite{Van}.

Digital and online learning have become an important part of modern education including in the field of engineering. Information technologies play a vital role in delivering digital learning to students. Colleges and universities have made significant investments to improve their information and communication technology (ICT) capacities \cite{Hernandez}. 

In response to the needs of Industry 4.0, some universities have adopted the framework of Education 4.0 \cite{Miranda, Ramirez}. The new education framework consists of four main components: i) competencies, ii) learning methods, iii) ICT, and iv) infrastructure. Students competencies are based on technological knowledge and skills for successful workplace performance, while the learning methods are based on problem solving and challenge-based learning. In particular, active and project-based learning  plays an important role in Education 4.0 \cite{Chen, Hernandez2}. Other innovative approaches to learning such as virtual-reality based engineering education can help improve the learning process related to Industry 4.0 \cite{Salah}.

In addition to technological progress, socio-cultural shifts must be taken into account in revising engineering curriculum. The new generation of students has its unique worldview which needs to be considered by the educators. In particular, the new generation is significantly affected by mobile devices and digital media. Educational content must be tailored to the new student preferences to achieve effective learning outcomes \cite{Moore}. Innovative approaches such gamification may help improve the learning process \cite{Mauricio, Ortiz}. 

Many universities have also introduced nontechnical updates to their engineering curriculum. The most significant nontechnical update has been the introduction of entrepreneurship courses and experiences for students. A lot of attention has recently been given to equipping students with entrepreneurial skills. Students learn about entrepreneurship in their courses as well as through university incubators.

\section{University mathematics curriculum}
The mathematics curriculum changed very little in the current century. It remains a largely analytic domain, where solutions are mainly obtained manually. The current mathematics curriculum emphasizes theory over practical approaches. For instance, when finding the extreme values of a function, derivative-based approach is preferred over the gradient decent. There are two key reasons for why analytical approaches are favored over numerical methods. First, analytical solutions are reliable and elegant. An analytical solution is guaranteed to be exact. Second, mathematics courses are usually taught by pure mathematicians who have an inherent preference for analytical solutions. Pure mathematics which is based on theorem proving is not amenable to numerical methods. 

Despite the popularity of analytical approaches to problem solving in mathematics, there has been a growing push to integrate computer algebra systems as part of the learning process. Computer algebra systems such as Matlab and Mathematica are now routinely used in many mathematics courses. 
%Ho Hon
The study by Cretchley et al. \cite{cretchley} found that engineering students were positive about the use of technology as a learning tool in mathematics courses. The increased use of technology in class helped improve student focus and interest in lectures. Student evaluations also indicated that they had a greater level of enjoyment towards the lectures due to the use of technology. It is noteworthy that students chose not to rely too heavily on technology during the examinations despite the freedom to do so. The students found it extremely important to be competent with analytical mathematical skills as opposed to purely computational skills. Some revealed that they learn the subjects equally well without the help of scientific packages, although the perception towards the use of computer is in general highly positive. Almost all students responded positively to Matlab as an effective tool for computation and graphing. Many used Matlab for non-examination purposes. For example, they utilized it to check their handwritten mathematical steps in assignments and practice problems; and others used it for exploration beyond the standard syllabus and curriculum.

The influence of computer technology on students' academic performance and learning experience has been investigated by several authors. Abdul Majid et al.\cite{majid0, majid} used Matlab as an aid to teach calculus to engineering students. The software package was used for various course learning outcomes such as graphical display of mathematical functions, exploration, identifying and predicting structural patterns in evaluating a series of complex indefinite integrals, and numerical approximations in applied mathematics. The study showed a positive impact on students' academic performance in the final examinations. The study concluded that the integration of scientific packages into engineering mathematics courses could be effective under certain conditions. Similarly, other studies \cite{Puhak, Strayhorn} also found a positive impact from the use of scientific software packages on students' motivation in learning mathematics. 

In a separate study by Brake \cite{Brake}, the authors investigated the use of Matlab in engineering mathematics courses to increase student confidence level and mathematical abilities. Matlab was used to solve concrete engineering problems which require a deep understanding of underlying mathematical principles. The study found generally positive student response to the use of software in their mathematics courses. However, the results of the study must be considered carefully given the small sample size of the subjects.

Although the majority of the studies were based on the use of Matlab, several other studies  considered alternative mathematics software packages.
The study by Kilicman et al. \cite{Kilicman} focused on the use of Maple to help students understand both the theoretical and computational aspects of linear algebra for engineering students. In particular, it was shown that the use of Maple facilitates the understanding of computational aspects of eigenvalues and eigenvectors. It allows students more time to focus on the theoretical aspects and the underlying mathematical principles. 

In a recent study by Mezhennaya and Pugachev \cite{Mez}, the authors compared engineering students' perceptions regarding several mathematical software: Matlab, Mathematica and Excel. 
The study found that all the scientific packages considered can be used in education, under the condition that the policies for software usage are carefully implemented. The study found that many students lack hands-on experience on how to use the software. The students particularly struggled with Matlab and Mathematica finding them non user friendly. The study concluded that additional classes are required to prepare students to use software in their courses.

\section{Numerical computing in mathematics curriculum}
Mathematics lies at the foundation of science and engineering. The importance of mathematics courses in engineering education cannot be underestimated. These courses equip students with the fundamental skills and knowledge to study the more specialized engineering courses. Thus, student success in engineering studies depends directly on the mathematics and sciences courses. Given the significance of the mathematics courses in the engineering curriculum, it is paramount to ensure their currency with respect to the Industry 4.0. 

The technological advances over the last decade have created demand for more computationally proficient experts. To meet this demand, numerical computing must become a core part of engineering studies. Mathematics courses offer a natural and convenient avenue for introducing numerical computing to engineering students. There are two main factors that make mathematics courses particularly amenable to numerical computing. First, in many cases mathematical problems have numerical solutions. For instance, finding the root of a polynomial or the minimum value of a function can be done numerically. Therefore, it is both logical and appropriate to apply numerical computing to mathematical problems.
Second, mathematics courses are usually taken at the beginning of the study plan. Thus, students become acquainted with numerical computing at an early stage. The computing and programming skills acquired in this manner will have a positive effect in the more advanced, downstream engineering courses.

The key idea for the proposed curriculum update is the addition of computing tutorials (labs) to mathematics courses. In particular, we propose adding weekly computing tutorials (labs) related to the main lecture material. For instance, in the week in which students cover finding the extreme values of a function, there will be a computing tutorial where students learn and implement the gradient descent algorithm. The suggested length of each tutorial is 1 hour. It is enough time to implement most of the numerical algorithms at the undergraduate level. At the same time, 1 extra hour per week will not overburden the students.

The exact details of numerical computing content is left for individual universities and instructors. Depending on the syllabus and course learning outcomes, the numerical computing labs will be different for each university and instructor. Nevertheless, the general ideas will be broadly similar across different curricula. To illustrate the proposed numerical computing content, we will focus on the three main concept in calculus: limits, derivatives, and integrals.

\subsection{Limits}
Limit is a fundamental concept in calculus. Students are usually taught to calculate limits using analytical approaches. Although analytical approaches work well, there is no single universal rule for calculating limits.
On the other hand, in most cases, limits can be calculated numerically using essentially the same approach. To illustrate, suppose we want to calculate  $\lim_{x\to a^+} f(x)$. Then we can loop for $k=0 \text{ to } n$ and calculate $f(a+10^{-k})$. 
As $k$ increases, $a+10^{-k}$ approaches $a$, so $f(a+10^{-k})$ will, in most cases, approach the limit value.
We can deduce the limit by observing the values of $f(a+10^{-k})$ or determine that the limit does not exist if there is no pattern of convergence. The value of $n$ can be chosen manually or using a stopping criterion. For instance, the algorithm may continue to iterate until the difference between consecutive values of $f(a+10^{-k})$ is below a certain threshold. The value of the limit can also be deduced automatically based on the values of $f(a+10^{-k})$ using various heuristics. 

Another common limit problem is $\lim_{x\to\infty} f(x)$. In this case, we can loop for $k=0 \text{ to } n$ and calculate $f(10^{k})$. As $k$ increases, $10^{k}$ approaches $\infty$, so $f(10^{k})$ will, in most cases, approach the limit value.
Then the limit can be determined based on the values of $f(10^{k})$. Various extensions and customizations of this basic approach can be made. For instance, to avoid issues with periodic functions $f(10^{k}+\epsilon_k)$, where $\epsilon_k$ are randomly generated, can be used. Other values than $10^k$ can also be used as long as the sequence approaches infinity. A degree of automization can be introduced using different heuristics.

\subsection{Derivative}
Derivative is arguably the most important concept in calculus. There exist several rules such as the power rule, the product rule, the chain rule, and others to find the derivative of a function by hand. However, manual differentiation may be cumbersome when dealing with complex function. On the other hand, calculating the derivative at a point numerically is relatively straightforward.
To illustrate, suppose that we want to calculate $f'(a)$. Recall that 
\begin{equation}
f'(a) = \lim_{x\to a} \frac{f(a+h)-f(a)}{h}.
\end{equation} 
Therefore, to calculate $f'(a)$ numerically we use the same approach as with the limits. In particular, we can loop for $k=0 \text{ to } n$ and calculate $\frac{f(a+10^{-k})-f(a)}{10^{-k}}$. Then the limit, and by extension the derivative, can be deduced (approximated) based on the calculated values. The accuracy of the approximation depends in large part on the value of $n$.

One of the most important applications of the derivative is finding the extreme values of a function. Traditionally, this is done by first finding the critical points of the function and then applying the second derivative test. However, finding the critical points is not always possible, so numerical approaches can be used in such cases. The most popular numerical approach for finding the extreme values is based on the gradient descent (ascent) algorithm. In gradient descent, the optimal value of $x$ is iteratively updated based on the gradient. In particular, for $k=0 \text{ to } n$, the updated optimal value of $x$ is given by
\begin{equation}
{x} _{k+1}= {x} _{k}-\alpha\nabla f({x} _{k}),
\end{equation}
where $\nabla f(x)$ is the gradient and $\alpha$ is the step size. In the case of a single-variable function, the gradient equals simply to the derivative $\nabla f(x) = f'(x)$. The step size $\alpha$ can be either fixed or dynamic. While a large value of $\alpha$  accelerates the convergence at the beginning, it may hurt the convergence in the region near the optimal value.

There exist several extensions of the basic gradient descent algorithm. One such extension is  gradient descent with momentum which uses the second derivative to anticipate the location of the next optimal point and thus accelerates the convergence.

\subsection{Integrals}
Integration is an important concept in engineering mathematics. Although there exist a number of rules for finding the integral, it is significantly more challenging than differentiation. Moreover, in many cases, the indefinite integral does not even exist. Therefore, numerical approaches are particularly useful for integration. 

To illustrate the application of numerical integration, suppose that we want to calculate $\int_a^b f(x)\,dx$. There exist several numerical methods for calculating the integral. One simple method is based on the Riemann sums. The interval $[a, b]$ is divided into $n$ equal subintervals with endpoints $x_0, x_1, ..., x_n$. Let $\delta = x_{k+1}-x_k$ be the length of each subinterval. Then the right Riemann sum is defined as $\delta\sum_{k=1}^n f(x_k)$. The Riemann sum provides an approximation of the integral. Indeed, 
\begin{equation}
\int_a^b f(x)\,dx = \lim_{n\to\infty} \delta\sum_{k=1}^n f(x_k).
\end{equation}
The Riemann sums can be quickly calculated on a computer providing a simple, yet effective approach to calculating integrals numerically. Other popular integral approximation methods include the trapezoid rule and the Simpson's rule.

\subsection{Additional considerations}

The above discussion about numerical methods for calculating limits, derivatives, and integrals is easily extended to multivariate calculus. For instance, to find the partial derivative $f_x(a, b)$, we can loop for $k=0 \text{ to } n$ and calculate $\frac{f(a+10^{-k}, b)-f(a, b)}{10^{-k}}$. Many problems related to sequences and series can similarly be solved using numerical techniques. In particular, the convergence of a series can be deduced from its partial sums. By calculating the partial sums on the computer and observing the results, we can intuit the nature of the series.

Vectorization is an important aspect of numerical computing. Since the modern computer chips are optimized for matrix multiplication, it is more efficient to employ vector operations. In particular, some algorithms based on for-loops can be converted into vector operations resulting in higher efficiency and speed. For instance, the Riemann sum can be calculated with a single vector operation:
\begin{equation}
\delta\sum_{k=1}^n f(x_k) =\delta S \big(f(\mathbf x)\big),
\end{equation}
where $\mathbf x = [x_1, ..., x_n]$ is the vector of endpoints, $f(\mathbf x)$ is a vectorized function operation, and $S$ is the vector function which returns the sum of all the coordinates.
Similarly, limit calculations can be vectorized and made more efficient. Vectorization is also useful in multi-variate calculus, where operations can be performed on a vector of variables.

The choice of the programming language for numerical computing requires careful consideration. There are several suitable candidates for this purpose including Python, Java, C++, Matlab, and others. Based on our experience with different programming languages, we recommend the use of Python. Python is currently the most popular programming language on the planet. It has a simple and intuitive syntax making it easy to learn and apply. Python has libraries to fit any purpose including  an extensive collection of libraries related to numerical computing. The basic Python libraries related to computing are NumPy, SciPy, and SymPy. More advanced packages such as OR-Tools are also available for optimization tasks.
Since Python is a universal programming language, it can be used for almost any task. Thus, students who learn Python in their mathematics courses can employ it in their other courses. In addition, numerical computing implemented in Python can be connected to other applications.

\section{Conclusion}
Although the classical approach to teaching mathematics is still relevant for certain student cohorts, it is outdated for engineering students. Modern engineering is increasingly reliant on computing \cite{Kamalov1, Thabtah}. Therefore, universities must equip the student with appropriate computing skills. In particular, mathematics courses must be revised to include numerical computing content.

Given the efficiency of computer-based calculations, numerical computing provides a convenient approach to problem solving in engineering mathematics.
It can be integrated into the existing curriculum with little hassle and cost. 
In this paper, we proposed a framework for integrating numerical computing into the existing mathematics curriculum. We demonstrated how numerical approaches can be used some of the most common problems encountered in calculus. The proposed framework can be customized by individual universities to fit their special needs.

\end{document}